\documentclass [12pt]{article}
\usepackage{amsfonts}
\usepackage{amsfonts}
\usepackage{amsfonts}
\usepackage{bbm}
\usepackage{amsfonts}
\usepackage{amsthm}
\usepackage{amsmath}
\usepackage{amstext}
\usepackage{amssymb}
\usepackage{mathrsfs}
\usepackage{epsf}              
\usepackage{graphicx}          
\usepackage{fancybox}          
\usepackage{color}             
\usepackage{fancyhdr}
\usepackage[hang,footnotesize]{caption2}  

\def\mathcal{\mathscr}
\newfont{\aaa}{cmb10 at 18pt}
\newfont{\bbb}{cmb10 at 10pt}

\numberwithin{equation}{section}
\newtheorem{thm}{Theorem}[section]
\newtheorem{defi}[thm]{Definition}
\newtheorem{cor}[thm]{Corollary}

\newtheorem{pro}[thm]{Proposition}

\newtheorem{rem}[thm]{Remark}

\def\pr{\noindent {\it Proof.}}

\def\QED{\hfill$\Box$}

\def\cl{\centerline}
\def\a{\alpha}

\def\la{\lambda}

\def\t{\tilde}
\def\B{{\mathcal B}}
\def\Z{\mathbb Z}
\def\CC{{\mathbb C}}
\def\pa{\partial}
\def\c{\circ}
\def\g{\gamma}
\def\Z{\mathbb{Z}}

\begin{document}

\setcounter{page}{1}
\qquad\\[8mm]

\cl{\large{\bf A Lie conformal algebra of Block type}}
 \vskip12pt
 \cl{Lamei Yuan}
 \cl{\small{Academy of Fundamental and Interdisciplinary
 Science,}}
\cl{\small{Harbin Institute of Technology, Harbin 150080, China}}  \cl{\small lmyuan@hit.edu.cn
 }\vskip6pt
{\small\parskip .005 truein \baselineskip 3pt \lineskip 3pt}

\noindent{\bf Abstract:} The aim of this paper is to study a Lie conformal algebra of Block type. In this paper, conformal derivation, conformal module of rank 1 and low-dimensional comohology of the Lie conformal algebra of Block type are studied. Also, the vertex Poisson algebra structure associated with the Lie conformal algebra of Block type is constructed.


\vskip5pt
\noindent{\bf MR(2000) Subject Classification:~}17B65, 17B69

\section{Introduction}

\noindent
The notion of Lie conformal algebra, introduced by
Kac \cite{Kac1}, encode an axiomatic description of the operator
product expansions of chiral fields
in conformal field theory. It is a powerful tool for the study of infinite-dimensional Lie
(super)algebras, associative algebras and their representations. Lie conformal algebras have been extensively studied, including the classification problem \cite{DK,FK}, cohomology theory \cite{BKV,Su1} and representation theory \cite{CK}.

The Lie conformal algebras are closely related to vertex algebras. Primc \cite{P} introduced and studied a notion of vertex Lie algebra, which is a special case of a more general
notion of local vertex Lie algebra \cite{DLM}. As it was explained in \cite{Li1}, the notion of Lie conformal algebra and the notion of vertex Lie algebra are equivalent. In this paper, we shall use Lie conformal algebra and vertex Lie algebra synonymously.

With the notion of vertex Lie algebra, one arrives at the notion of vertex
Poisson algebra, which is a combination of a differential algebra structure and a vertex Lie algebra structure, satisfying a natural compatibility condition. The symmetric algebra of
a vertex Lie algebra is naturally a vertex Poisson algebra \cite{FB}. A general construction theorem of vertex Poisson algebras was given in \cite{Li1}. Applications of vertex Poisson algebras to the
theory of integrable systems were studied in \cite{BDSK}.

In the present paper, we study a nonsimple Lie conformal algebra of infinite rank, which is endowed with a $\CC[\pa]$-basis $\{J_i|i\in\Z^+\}$, such that
\begin{eqnarray}\label{bracket}
[J_{i\,\la} J_j]=((i+1)\pa+(i+j+2)\la)J_{i+j}, \ \mbox{for}\
i, j\in \Z^+.
\end{eqnarray}
The corresponding formal distribution Lie algebra is a Block type Lie algebra, which is the associated graded Lie algebra of the filtered Lie algebra $W_{1+\infty}$ \cite{Su2, Su3, SXX1,SXX2,WT}. Thus we call this Lie conformal algebra a {\it Lie conformal algebra of Block type} and denote it by $\B$ in this paper. It is a conformal subalgebra of ${\rm gr}\, gc_1$
studied in \cite{SY1}. In addition, it contains the Virasoro conformal algebra $\rm Vir=\CC[\pa]J_0$ with $[J_{0\,\la} J_0]=(\pa+2\la)J_0$ as a subalgebra.

The paper is organized as follows. In Section 2, we recall the notions of Lie conformal algebra and vertex Lie algebra. In Section 3, we study conformal derivations of the Lie conformal algebra of Block type $\B$. In Section 4, we recall the notions of conformal module and comohology of Lie conformal algebras. Then we study conformal module of rank 1 and low-dimensional comohology of $\B$ with coefficients in $\B$-modules. In Section 5, we equip a vertex Lie algebra structure $(Y_-,\pa)$ with $\B$ and establish an association of a vertex Poisson algebra structure to the vertex Lie algebra $(\B, Y_-,\pa)$.

\section{Preliminaries}

\noindent
Throughout this paper, all vector spaces and tensor products are over the complex field $\mathbb{C}$.  We use notations $\Z$ for the set of integers and $\Z^+$ for the set of nonnegative integers.

\begin{defi}\rm A Lie conformal algebra $\mathcal {R}$ is a $\CC[\partial ]$-module with a $\CC$-bilinear map,
\begin{eqnarray*}
\mathcal {R}\otimes \mathcal {R}\rightarrow \CC[\lambda]\otimes \mathcal {R}, \ \ \ a\otimes b \mapsto [a_\lambda b],
\end{eqnarray*}
called the $\lambda$-bracket, and
satisfying the following axioms ($a, b, c\in \mathcal {R}$),
\begin{eqnarray}
\mbox{(conformal\  sesquilinearity)} \ \ [\partial a_\lambda b]&=&-\lambda[a_\lambda b],\ [ a_\lambda \partial b]=(\partial+\lambda)[a_\lambda b],\label{L1}\\
 \mbox{(skew-symmetry)}\ \ \ \  {[a_\lambda b]} &=& -[b_{-\lambda-\partial}a], \label{L2}\\
\mbox{(Jacobi \ identity)}\ \ {[a_\lambda[b_\mu c]]}&=&[[a_\lambda b]_{\lambda+\mu
}c]+[b_\mu[a_\lambda c]].\label{L3}
\end{eqnarray}
\end{defi}
If we consider the expansion
\begin{eqnarray}\label{j-product}[a_{\lambda}b]=\mbox{$\sum\limits_{j\in\Z^+
}$}\frac{\lambda^j}{j!}(a_{(j)}b),\end{eqnarray}
the coefficients of $\frac{\lambda^j}{j!}$ are called the {\bf $j$-product} satisfying ${a_{(n)}b}=0$ for $n$ sufficiently large,
and the axioms \eqref{L1}--\eqref{L3} can be written in terms of them as follows:
\begin{eqnarray}
{\partial a_{(n)} b}&=&-na_{(n-1)} b, \ \ {a_{(n)} \partial b}=\partial(a_{(n)}b)+n
a_{(n-1)}b, \label{jHL1}\\
{{a_{(n)}b}} &=& -\mbox{$\sum\limits _{i\in\Z^+}$}(-1)^{n+i}\frac{1}{i!}\pa^i b_{(n+i)}a,\label{jHL2}\\
{a_{(m)}b_{(n)} c}&=&{b_{(n)}a_{(m)}
c}+\mbox{$\sum\limits_{i=0}^{m}$}\mbox{${m\choose
i}$}(a_{(i)}b)_{(m+n-i)}c.\label{jjHL3}
\end{eqnarray}

In terms of generating functions, Mirko Primc in \cite{P} presented an equivalent
definition of a Lie conformal algebra under the name of vertex Lie algebra (see also \cite{Li1}). Let $V$ be any vector space.
Following \cite{P}, for a formal series
\begin{eqnarray*}
f(x_1,\cdots,x_n)=\mbox{$\sum\limits_{m_1,\cdots,m_n\in\Z}$}u(m_1,\cdots,m_n)x_1^{-m_1-1}\cdots x_n^{-m_n-1}\in V[[x_1^{\pm 1},\cdots, x_n^{\pm 1}]],
\end{eqnarray*}
we set
\begin{eqnarray}\label{sing}
{\rm Sing} f(x_1,\cdots,x_n)=\mbox{$\sum\limits_{m_1,\cdots,m_n\in\Z^+}$}u(m_1,\cdots,m_n)x_1^{-m_1-1}\cdots x_n^{-m_n-1}.
\end{eqnarray}
Clearly, for $1\leq i\leq n$,
\begin{eqnarray}\label{sing1}
\frac {\pa}{\pa x_i}{\rm Sing}f(x_1,\cdots,x_n)={\rm Sing}  \frac {\pa}{\pa x_i} f(x_1,\cdots,x_n).
\end{eqnarray}

\begin{defi}\label{d2}{\rm A vertex Lie algebra is a vector space
$A$ equipped with a linear operator $\partial$ called the derivation and a
linear map
\begin{eqnarray*}
Y_-(\cdot, z): A \rightarrow z^{-1}({\rm End\,}
A)[[z^{-1}]],\ \ \ a\mapsto Y_-(a,z)=\mbox{$\sum_{n\geq
0}$}a_{(n)}z^{-n-1},
\end{eqnarray*}
satisfying the following conditions for $a,b \in A $, $n
\in\Z^+$:
\begin{eqnarray}
&&{a_{(n)}b}=0 \ \ \mbox{for} \ n \ \mbox{sufficiently large},\label{vHL0+}\\
&&{[\pa, Y_-(a,z)]}={Y_-(\pa a,z)}=\frac{d}{dz}Y_-(a,z),\label{vHL1+}\\
&&{Y_-(a,z)b}= {\rm Sing} \big(e^{z\pa}Y_-(b,-z)a\big),\label{vHL2+}
\end{eqnarray}
and the half Jacobi identity holds:
\begin{eqnarray}\label{vHL3+}
&&{\rm Sing}\Big(z_0^{-1}\delta(\frac{z_1-z_2}{z_0})Y_-(a,z_1)Y_-(b,z_2)-z_0^{-1}\delta(\frac{z_2-z_1}{-z_0})Y_-(b,z_1)Y_-(a,z_2)\Big) \nonumber\\
& &\ \ \ \ \ \ \ \ \ \ \hspace{2cm} ={\rm Sing} \Big(z_2^{-1}\delta(\frac{z_1-z_0}{z_2})Y_-(Y_-(a,z_0)b,z_2)\Big).
\end{eqnarray}}
\end{defi}
 Relation \eqref{vHL2+} is called the {\it half skew-symmetry}.
It was shown in \cite{P} that the half Jacobi identity \eqref{vHL3+} amounts to the following {\it half commutator formula}:
\begin{eqnarray}\label{vHL55}
Y_-(a,z_1)Y_-(b,z_2)-Y_-(b,z_1)Y_-(a,z_2)
={\rm Sing}\Big(\mbox{$\sum_{i\in\Z^+}$}(z_1-z_2)^{-i-1}Y_-(a_{(i)}b,z_{2})\Big).
\end{eqnarray}
As it was explained in \cite[Remark 2.6]{Li1}, the notion of Lie conformal algebra is equivalent to the notion of vertex Lie algebra. We often denote a vertex Lie algebra by $(A, Y_-,\pa)$ and refer to $(Y_-,\pa)$ as the vertex Lie algebra structure. A vertex Lie algebra $(A, Y_-,\pa)$ is said to be {\it free}, if $A$ is a free $\CC[\pa]$-module over a vector space $V$, namely, $A=\CC[\pa]V\cong \CC[\pa]\otimes_{\CC} V.$

\section{Conformal derivation}

\vskip6pt
\noindent

Let $\mathcal {C}$ denote the ring $\mathbb{C}[\partial]$
of polynomials in the indeterminate $\partial$.
\begin{defi}\rm Let $V$ and $W$ be two $\mathcal {C}$-modules.
A linear map $\phi:V\rightarrow \mathcal
{C}[\lambda]\otimes_{\mathcal {C}}W$, denoted by $\phi_\lambda: V
\rightarrow W $, is called a conformal linear map, if
\begin{eqnarray}
\phi_{\lambda}(\partial v)=(\partial+\lambda)(\phi_\lambda v),\ \mbox{for}\ v\in V.\end{eqnarray}
\end{defi}
The space of conformal
linear maps between $\mathcal{C}$-modules $V$ and $W$ is denoted by
${\rm Chom}(V,W)$ and it can be made into an $\mathcal {C}$-module
via
\begin{eqnarray*}
(\pa\phi)_\la v=-\la\phi_\la v, \ \mbox{for} \ v\in V.
\end{eqnarray*}
\begin{defi}\rm Let $\mathcal {R}$ be a Lie conformal algebra.
A conformal linear map $d_\lambda:\mathcal {R}\rightarrow \mathcal
{R}$ is called a conformal derivation of $\mathcal {R}$ if
\begin{eqnarray}\label{der}
d_\lambda[a_\mu b]=[(d_\lambda a)_{\lambda+\mu}b]+[a_\mu(d_\lambda b)],\  \mbox{for}\ a, b\in\mathcal {R}.
\end{eqnarray}
\end{defi}
The space of all conformal derivations of $\mathcal {R}$ is denoted
by ${\rm CDer}(\mathcal {R})$. For any $a\in\mathcal {R}$, one can
define a conformal derivation $({\rm ad}\,a)_\lambda:\mathcal
{R}\rightarrow \mathcal {R}$ by $({\rm ad}\,a)_\lambda
b=[a_\lambda b]$ for $b\in\mathcal {R}$. Such conformal derivation is called   {\it inner}. Denote by ${\rm CInn}(\mathcal {R})$ the space of all conformal inner derivations of $\mathcal {R}$.

\begin{pro}
Every conformal derivation of the Lie conformal algebra $\B$ is
inner.
\end{pro}
\noindent{\it Proof.~}Let $d$ be any conformal derivation of $\B$. Denote $L=J_0$. Assume that there exists a finite subset $I=\{i_1,\cdots,i_n\}\subseteq \Z^+$ such that
$d_\la L=\mbox{$\sum_{j=1}^{n}$}f_{i_j}(\pa,\la)J_{i_j}$, where $f_{i_j}(\pa,\la)\in\mathbb{C}[\pa,\la]$.
Condition \eqref{der} requires $d_\la([L_\mu L]) =
[L_\mu (d_\la L)]+[(d_\la L)_{\la+\mu} L].$ This is equivalent to
\begin{eqnarray}
&&(\pa+\la+2\mu)\mbox{$\sum\limits_{j=1}^{n}$}f_{i_j}(\pa,\la)J_{i_j}-\mbox{$\sum\limits_{j=1}^{n}$}(\pa+(i_j+2)\mu)f_{i_{j}}(\pa+\mu,\la)J_{i_j}\nonumber \\&&\ \ \ \ =\mbox{$\sum\limits_{j=1}^{n}$}((i_j+1)\pa+(i_j+2)(\la+\mu))f_{i_j}(-\la-\mu,\la)J_{i_j}.
\end{eqnarray}
For each $j$,
\begin{eqnarray}\label{der1}
&&(\pa+\la+2\mu)f_{i_j}(\pa,\la)-(\pa+(i+2)\mu)f_{i_j}(\pa+\mu,\la)\nonumber\\&&\ \ \ \ =((i_j+1)\pa+(i_j+2)(\la+\mu))f_{i_j}(-\la-\mu,\la).
\end{eqnarray}
Write $f_{i_j}(\lambda,\partial)=\sum_{k=0}^{m}a_{i_j,k}(\lambda)\partial^k$
with $a_{i_j,m}(\la)\neq 0$. Then, assuming $m > 1$, if we equate terms
of degree $m$ in $\partial$, we have
$(\la-i-m\mu)a_{i_j,m}(\la)=0$ and thus $a_{i_j,m}(\la)=0$. This
contradicts $a_{i_j,m}(\la)\neq 0$. Thus ${\rm deg}_{\pa} f_{i_j}(\lambda,\partial)\leq 1$, and  $f_{i_j}(\lambda,\partial)=a_{i_j,0}(\la)+a_{i_j,1}(\la)\pa$. Substituting it into \eqref{der1} gives
$a_{i_j,0}(\la)=\frac{i_j+2}{i_j+1}\la a_{i_j,1}(\la)$. Therefore,
\begin{eqnarray*}
d_\la
L=\mbox{$\sum\limits_{j=1}^{n}$}\frac{a_{i_j,1}(\la)}{i_j+1}((i_j+1)\pa+(i_j+2)\la)J_{i_j}.
\end{eqnarray*}
Replacing $d_\la$ by $d_\la-({\rm ad}\,h)_\lambda$ with $h=\mbox{$\sum_{j=1}^{n}$}\frac{a_{i_j,1}(-\pa)}{i_j+1}J_{i_j}$, we get $d_\la(L)=0$.

For $k>0$, assume that $d_\la J_k=\sum_{i=1}^{l}f_{k_i}(\pa,\la)J_{k_i}.$
Applying $d_\la$ to $[L_\mu J_k]=(\pa+(k+2)\mu)J_{k}$ and using
$d_\la(L)=0$, we obtain
\begin{eqnarray}
(\pa+\la+(k+2)\mu)\mbox{$\sum\limits_{i=1}^{l}$}f_{k_i}(\pa,\la)J_{k_i}
=\mbox{$\sum\limits_{i=1}^{l}$}(\pa+(k_i+2)\mu)f_{k_i}(\pa+\mu,\la)J_{k_i},
\end{eqnarray}
and thus
\begin{eqnarray}
(\pa+\la+(k+2)\mu)f_{k_i}(\pa,\la)
=(\pa+(k_i+2)\mu)f_{k_i}(\pa+\mu,\la), \ \ \mbox{for}\ 1\leq i\leq
l.
\end{eqnarray}
Comparing the highest degree of $\la$ gives $f_{k_i}(\pa,\la)=0$ for $ 1\leq i\leq l$. Hence $d_\la(J_k)=0$ for $k>0$. This concludes the proof. \QED

\section{Cohomology}

\begin{defi}\rm A module $M$ over a Lie conformal algebra $\mathcal {R}$ is a
$\mathbb{C}[\partial]$-module endowed with a bilinear map
\begin{eqnarray*}
\mathcal {R}\otimes M\rightarrow
M[[\lambda]],\ \ a\otimes v\mapsto a_\lambda v
\end{eqnarray*}
such that ($a,b\in\mathcal {R}$, $v\in M$)
\begin{eqnarray}
&&a_\lambda(b_\mu v)-b_\mu(a_\lambda v)=[a_\lambda b]_{\lambda+\mu}v,\\
&&(\partial a)_\lambda v=-\lambda a_\lambda v,\ \ \
a_\lambda(\partial v)=(\partial+\lambda)a_\lambda v,
\end{eqnarray}
If $a_\lambda v\in
M[\lambda]$ for all $a \in \mathcal {R}$, $v \in M$, then the
$\mathcal {R}$-module $M$ is said to be conformal. If $M$ is
finitely generated as $\mathbb{C}[\partial]$-module, then $M$ is simply
called finite.
\end{defi}

Since we only consider conformal modules, we will simply shorten the term
``conformal module" to ``module". The one-dimensional vector space $\CC$ can be viewed as a module (called the {\it trivial module}) over any conformal algebra $\mathcal {R}$
with both the action of $\pa$ and the action of $\mathcal {R}$ being
zero. In addition, for a fixed nonzero complex constant $a$, there
is a natural $\CC[\pa]$-module $\CC_a$, which is the one-dimensional
vector space $\CC$ such that $\pa v=a v$ for $v\in\CC_a$. Then
$\CC_a$ becomes an $\mathcal{R}$-module on which all elements of
$\mathcal{R}$ act by zero.

For the Virasoro conformal algebra Vir, it was proved in \cite{CK} that all free nontrivial
Vir-modules of rank $1$ are the
following ones $(\Delta, \alpha\in \mathbb{C})$:
\begin{eqnarray}\label{Vir-m}
M_{\Delta,\alpha}=\mathbb{C}[\partial]v,\ \ L_\lambda
v=(\partial+\alpha+\Delta \lambda)v.
\end{eqnarray}
The module $M_{\Delta,\alpha}$ is irreducible if and only if
$\Delta\neq 0$, the module $M_{0,\alpha}$ contains a unique
nontrivial submodule $(\partial +\alpha)M_{0,\alpha}$ isomorphic to
$M_{1,\alpha}$, and the modules
$M_{\Delta,\alpha}$ with $\Delta\neq 0$ exhaust all finite
irreducible nontrivial conformal Vir-modules.
\begin{pro}\label{theo} All free nontrivial $\B$-modules of rank $1$ are as follows $(\Delta, \alpha\in \mathbb{C})$:
\begin{eqnarray*}
M_{\Delta,\alpha}=\mathbb{C}[\partial]v,\ J_{0_\lambda} v=(\partial+\alpha+\Delta
\lambda)v, \ J_i\,{_\lambda} v=0,\ \mbox{for}\ i> 0.\end{eqnarray*}
\end{pro}
\noindent{\it Proof.} By \eqref{Vir-m}, $J_{0\,\lambda} v=(\partial+\alpha+\Delta
\lambda)v$ for some $\Delta, \alpha\in \mathbb{C}$. By \cite[Lemma 5.1]{SY1}, we can suppose
that $k$ is the smallest nonnegative integer
such that $J_{k}\,_{\la}v\neq0$, $J_{k+1}\,_{\la} v=0$. Assume $k>0$ and write $J_{k}\,_{\la} v=g(\la,\pa)v$, for some
$g(\la,\pa)\in\CC[\la,\pa]$. Since $[J_{k}\,_{\la} J_{k}]_{\la+\mu}v=0$,
\begin{eqnarray}
g(\la,\pa)g(\mu,\la+\pa)=g(\mu,\pa)g(\la,\mu+\pa).
\end{eqnarray}
This implies ${\rm deg}_\la g(\la,\pa)+{\rm deg}_\pa
g(\la,\pa)={\rm deg}_\la g(\la,\pa)$. Thus ${\rm deg}_\pa g(\la,\pa)=0$. Then we have
$g(\la,\pa)=g(\la)$ for some $g(\la)\in\CC[\la]$. The fact that
$[J_{0\,\la} J_k]_{\la+\mu}v=((k+1)\la-\mu) J_k{_{\la+\mu}}v$ yields
\begin{eqnarray*}
((k+1)\la-\mu)
g(\la+\mu)v=(\pa+\a+\Delta\la)g(\mu)v-(\pa+\mu+\a+\Delta\la)g(\mu)v=
-\mu g(\mu)v,
\end{eqnarray*}
which gives $g(\mu)=0$. Hence, $J_{k}\,_{\la} v=0$, a
contradiction. Thus $k=0$ and $J_{1}\,_{\la} v=0$.
It follows immediately that $J_{i}\,_{\la} v=0$ for all $i\geq 1$.\QED\vskip5pt

In the following we study cohomology of the Lie conformal algebra $\B$
with coefficients in $\B$-modules $\CC$, $\CC_a$ and
$M_{\Delta,\alpha}$, respectively. For completeness, we shall present the definition of cohomology of Lie conformal algebras given in \cite{BKV}.

\begin{defi}\label{cochain} \rm An $n$-cochain ($ n\in\Z^+$) of a Lie conformal algebra $\mathcal{R}$ with coefficients in an
$\mathcal{R}$-module $M$ is a $\CC$-linear map
\begin{eqnarray*}
\gamma:\mathcal{R}^{\otimes n}\rightarrow M[\la_1,\cdots,\la_n],\ \
\ a_1\otimes\cdots \otimes a_n \mapsto
\g_{\la_1,\cdots,\la_n}(a_1,\cdots,a_n)
\end{eqnarray*}
satisfying

(1) $\g_{\la_1,\cdots,\la_n}(a_1,\cdots,\pa a_i,\cdots,
a_n)=-\la_i\g_{\la_1,\cdots,\la_n}(a_1,\cdots, a_n)$ \ \mbox{(conformal antilinearity)},

(2) $\g$ is skew-symmetric with respect to simultaneous permutations of $a_i$'s and $\la_i$'s, namely,
\begin{eqnarray}
&& \ \ \ \ \ \ \ \ \ \g_{\la_1,\cdots,\la_{i-1},\la_{i+1},\la_i,\la_{i+2},\cdots,\la_n}(a_1,\cdots,a_{i-1},a_{i+1},a_i,\a_{i+2},\cdots,
a_n)\nonumber\\ &&\ \ \ \ \ \ \ \ \ \ \ \ \ \ =-\g_{\la_1,\cdots,\la_{i-1},\la_i,\la_{i+1},\la_{i+2},\cdots,\la_n}(a_1,\cdots,a_{i-1},a_i,a_{i+1},a_{i+2},\cdots,
a_n).\label{skey-sym}
\end{eqnarray}
\end{defi}

As usual, let $\mathcal{R}^{\otimes 0}= \CC$, so that a $0$-cochain
is an element of $M$. Denote by  ${\t C}^n(\mathcal {R},M)$ the set
of all $n$-cochains. The differential $d$ of an $n$-cochain $\gamma$ is
defined by
\begin{eqnarray}\label{d}
&&(d\g)_{\la_1,\cdots,\la_{n+1}}(a_1,\cdots,a_{n+1})\nonumber\\ &&=\mbox{$\sum\limits_{i=1}^{n+1}$}(-1)^{i+1}a_{i\,{\la_i}}\g_{\la_1,\cdots,\hat{\la_i},\cdots,\la_{n+1}}(a_1,\cdots,\hat{a_i},\cdots,a_{n+1})\nonumber\\
&&+\mbox{$\sum\limits_{i,j=1; i<j}^{n+1}$}(-1)^{i+j}\g_{\la_i+\la_j,\la_1,\cdots,\hat{\la_i},\cdots,\hat{\la_j},\cdots,\la_{n+1}}
\big([a_{i\,{\la_i}}a_j],a_1,\cdots,\hat{a_i},\cdots,\hat{a_j},\cdots,a_{n+1}\big),
\end{eqnarray}
where $\g$ is extended linearly  over the polynomials in $\la_i$. In
particular, if $\g\in M$ is a $0$-cochain, then
$(d\g)_\la(a)=a_\la\g$.

It is proved in \cite{BKV} that the operator $d$ preserves the space of cochains and $d^2=0$. Thus the cochains of a Lie conformal algebra $\mathcal{R}$ with coefficients in $\mathcal{R}$-module $M$ form a complex, which is called the {\it basic complex} and will be denoted by
\begin{eqnarray*}
\t C^\bullet(\mathcal{R},M)=\mbox{$\bigoplus\limits_{n\in\Z^+}$}\t C^n(\mathcal{R},M).
\end{eqnarray*}
Moreover, define a $\CC[\pa]$-module structure on $\t
C^\bullet(\mathcal{R},M)$ by
\begin{eqnarray}\label{co-mod}
(\pa\g)_{\la_1,\cdots,\la_n}(a_1,\cdots,
a_n)=(\pa_M+\mbox{$\sum\limits_{i=1}^n$}\la_i)\g_{\la_1,\cdots,\la_n}(a_1,\cdots,
a_n),
\end{eqnarray}
where $\pa_M$ denotes the action of $\pa$ on $M$. Then $d\pa=\pa d$
and thus $\pa \t C^\bullet(\mathcal{R},M)\subset \t
C^\bullet(\mathcal{R},M)$ forms a subcomplex. The quotient
complex
\begin{eqnarray*}
C^\bullet(\mathcal{R},M)=\t C^\bullet(\mathcal{R},M)/\pa \t
C^\bullet(\mathcal{R},M)= \mbox{$\bigoplus\limits_{n\in\Z^+}$}
C^n(\mathcal{R},M)
\end{eqnarray*}
is called the {\it reduced complex}.
\begin{defi}\label{def11}\rm The basis cohomology $\t H^\bullet (\mathcal{R},M)$ of a Lie conformal algebra $\mathcal{R}$ with coefficients
 in $\mathcal{R}$-module $M$ is the
cohomology of the basis complex $\t C^\bullet(\mathcal{R},M)$ and
the (reduced) cohomology $H^\bullet (\mathcal{R},M)$ is the
cohomology of
the reduced complex $C^\bullet(\mathcal{R},M)$.
\end{defi}
\begin{rem}{\rm The basic cohomology $\t H^\bullet (\mathcal{R},M)$ is naturally a $\CC[\pa]$-module, whereas the reduced cohomology $H^\bullet (\mathcal{R},M)$ is a complex vector space.}\end{rem}
For a $q$-cochain $\gamma\in{\t C}^q(\mathcal {R},M)$, we call
$\gamma$ a {\it $q$-cocycle} if $d(\gamma)=0$; a {\it $q$-coboundary} or a
{\it trivial $q$-cocycle} if there is a $(q-1)$-cochain $\phi\in\t
C^{q-1}(\mathcal{R},M)$ such that $\gamma=d(\phi)$. Two cochains
$\gamma_1$ and $\gamma_2$ are called {\it equivalent} if $\gamma_1-\gamma_2$ is a
coboundary. Denote by $\t D^q(\mathcal{R},M)$ and $\t B^q(\mathcal{R},M)$ the spaces
of $q$-cocycles and $q$-boundaries, respectively. By Definition \ref{def11},
\begin{eqnarray*}
{\rm \t H}^q(\mathcal{R},M)=\t D^q(\mathcal{R},M)/\t B^q(\mathcal{R},M)=\{\mbox{equivalent classes of
$q$-cocycles}\}.
\end{eqnarray*}

The main results of this section are the following theorem.
\begin{thm}\label{thm-5}  For the Lie conformal algebra $\mathcal{B}$, the
following statements hold.\begin{itemize}\parskip-3pt
\item[\rm(1)] For the trivial module $\CC$, we have
\begin{eqnarray*}
{\rm dim\,\t H}^q(\mathcal{B},\CC)=\left\{
\begin{array}{ll}
1 &{\mbox if}\ q=0,\\
0 &{\mbox if}\ q=1, {\mbox or} \ 2,
\end{array}
\right.
\end{eqnarray*}
and
\begin{eqnarray*}
{\rm dim\, H}^q(\mathcal{B},\CC)=\left\{
\begin{array}{ll}
1 &{\mbox if}\ q=0, {\mbox or}\ 2, \\
0 &{\mbox if}\ q=1.
\end{array}
\right.
\end{eqnarray*}
\item[\rm(2)] If $a\neq 0$, then  ${\rm dim\,H}^\bullet(\B,\CC_a)=0$.
\item[\rm(3)] If $\a\neq 0$, then
$ {\rm dim\,H}^\bullet(\mathcal{B},M_{\Delta,\a})=0$.
\end{itemize}
\end{thm}
\noindent{\it Proof.~} (1) Since a $0$-cochain $\gamma$ is an
element of $\CC$ and $(d\gamma)_\la (a)=a_\la \gamma =0$ for
$a\in \mathcal{B}$, we have $\t D^0(\mathcal{B},\CC)=\t
C^0(\mathcal{B},\CC)=\CC$ and $\t B^0(\mathcal{B},\CC)=0$. Thus $\t
{\rm H}^0(\mathcal{B},\CC)=\t D^0(\mathcal{B},\CC)/\t B^0(\mathcal{B},\CC)=\CC$, and ${\rm H}^0(\mathcal{B},\CC)=\CC$ because $\pa \t C^0(\mathcal{B},\CC)=\pa\CC=0$.

Let $\gamma\in\t C^1(\mathcal{B},\CC)$ and
$d\gamma\in\pa\t C^2(\mathcal{B},\CC)$, namely, there is $\phi\in\t
C^2(\mathcal{B},\CC)$ such that $d(\gamma)=\partial \phi$. By \eqref{d}, \eqref{co-mod} and $\pa\CC=0$,
\begin{eqnarray}\label{5--1}
\g_{\la_1+\la_2}([a_{\la_1} b])=-(d\g)_{\la_1,\la_2}(a,b)=-(\pa
\phi)_{\la_1,\la_2}(a,b)=-(\la_1+\la_2)\phi_{\la_1,\la_2}(a,b), \ \ a,b\in\B,
\end{eqnarray}
By \eqref{bracket}, \eqref{5--1} and Definition \ref{cochain} (1),
\begin{eqnarray}
((i+1)\la_1-\la_2)\g_{\la_1+\la_2}(J_i)=-(\la_1+\la_2)\phi_{\la_1,\la_2}(J_0,J_i), \ i\geq0. \label{5**}
\end{eqnarray}
Setting $\la=\la_1+\la_2$ in \eqref{5**} gives
\begin{eqnarray*}
((i+1)\la-(i+2)\la_2)\g_{\la}(J_i)=-\la\phi_{\la_1,\la_2}(L,J_i),  \ i\geq0,\label{7-1}
\end{eqnarray*}
 which implies that $\g_{\la}(J_i)$ is divisible by $\la$. We can define a $1$-cochain $\g'\in\t C^1(\B,\CC)$ by
\begin{eqnarray}\label{7-12}
\g'_{\la}(J_i)=\la^{-1}\g_{\la}(J_i), \ \mbox{for}\ i\geq 0.
\end{eqnarray}
Since $\pa\CC=0$, $\g=\pa \g'\in\pa\t
C^1(\B,\CC)$. Hence ${\rm H}^1(\B,\CC)=0$. If $\g$ is a 1-cocycle, namely, $\phi=0$ in \eqref{5--1}, then \eqref{5**} gives $\g=0$. Thus ${\rm \t H}^1(\B,\CC)=0$.

Let $\psi\in\t D^2(\B,\CC)$ be a 2-cocycle. We have
\begin{eqnarray}
0&=&(d\psi)_{\la_1,\la_2,\la_3}(J_i,J_0,J_0)|_{\la_3=0}\nonumber\\
&=& -(\la_1-(i+1)\la_2)\psi_{\la_1+\la_2,\la_3}(J_i,J_0)|_{\la_3=0}
+(\la_1-(i+1)\la_3)\psi_{\la_1+\la_3,\la_2}(J_i,J_0)|_{\la_3=0}\nonumber\\
&&-(\la_2-\la_3)\psi_{\la_2+\la_3,\la_1}(J_0,J_i)|_{\la_3=0}\nonumber
\\&=&-((\la_1-(i+1)\la_2)\psi_{\la_1+\la_2,0}(J_i,J_0)+(\la_1+\la_2)\psi_{\la_1,\la_2}(J_i,J_0).\label{5--2}
\end{eqnarray}
Setting $\la=\la_1+\la_2$ in \eqref{5--2} gives $((\la-(i+2)\la_2)\psi_{\la,0}(J_i,J_0)=\la\psi_{\la_1,\la_2}(J_i,J_0)$. Thus $\psi_{\la,0}(J_i,J_0)$ is divisible by $\la$.
Define
a $1$-cochain $f$ by
\begin{eqnarray*}
f_{\la_1}(J_i)=\la_1^{-1}\psi_{\la_1,\la}(J_i,J_0)|_{\la=0}, \ \mbox{for}\ i\geq0.
\end{eqnarray*}
Set $\g=\psi+d f$, which is also a 2-cocycle. For all $i\geq0$,
\begin{eqnarray}
\g_{\la_1,\la}(J_i,J_0)|_{\la=0}=\psi_{\la_1,\la}(J_i,J_0)|_{\la=0}-\la_1f_{\la_1}(J_i)=0.\label{5++}
\end{eqnarray}
By \eqref{5++} and \eqref{5--2} with $\g$ in place of $\psi$, we have
$(\la_1+\la_2)\g_{\la_1,\la_2}(J_i,J_0)=0.$
Therefore $\g_{\la_1,\la_2}(J_i,J_0)=0=\g_{\la_1,\la_2}(J_i,J_0)$. With this,
\begin{eqnarray}
0&=&(d\g)_{\la_1,\la_2,\la}(J_0,J_i,J_k)|_{\la=0}\nonumber\\
&=&-\g_{\la_1+\la_2,\la}([J_{0\,{\la_1}}J_i],J_k)|_{\la=0}+\g_{\la_1+\la,\la_2}([J_{0\,{\la_1
}}J_k],J_i)|_{\la=0}-\g_{\la_2+\la,\la_1}([J_{i\,{\la_2}} J_k],J_0)|_{\la=0}\nonumber\\
&=&-((i+1)\la_1-\la_2)\g_{\la_1+\la_2,0}(J_i,J_k)-(k+1)\la_1\g_{\la_2,\la_1}(J_i,J_k).\label{7-2}
\end{eqnarray}
Setting $\la_1=0$ in \eqref{7-2} gives
$\g_{\la_2,0}(J_i,J_k)=0$ and thus
$\g_{\la_2,\la_1}(J_i,J_k)=0$. This proves $\g=0$. Hence ${\rm \t H}^2(\B,\CC)=0$.

It remains to compute ${\rm H}^2(\B,\CC)$. Following \cite{Su1}, we define a linear map $\sigma\gamma:\mathcal{B}^{\otimes q}\rightarrow \CC[\la_1,\cdots,\la_{q-1}]$ for $q\geq 2$ by
\begin{eqnarray}\label{sigma}
(\sigma\gamma)(\ a_1\otimes\cdots \otimes a_n)=\g_{\la_1,\cdots,\la_q}(a_1,\cdots,a_q)|_{\la_q=-\la_1-\cdots-\la_{q-1}}, \ \ a_1,\cdots,a_q\in\B.
\end{eqnarray}
 We define $\sigma\gamma=\g$ if $q=0$ and $\sigma\gamma(a_1)=\g_{\la_1}(a_1)|_{\la_1=0}$ if $q=1$. Set $C'^q(\mathcal{B},\CC)=\{\sigma\gamma|\gamma\in{\t C}^q(\mathcal {B},\CC)\}$.
Obviously, $\sigma:{\t C}^q(\mathcal {B},\CC)\rightarrow C'^q(\mathcal{B},\CC)$ is a surjective map. If $\g\in\pa {\t C}^q(\mathcal {B},\CC)=(\mbox{$\sum_{i=1}^q$}\la_i)\t
C^q(\B,\CC)$, then $\sigma\gamma=0$. That is, $\sigma$ factors to a map $\sigma:C^q(\mathcal {B},\CC)\rightarrow C'^q(\mathcal{B},\CC)$.

We claim that $\sigma:C^q(\mathcal {B},\CC)\rightarrow C'^q(\mathcal{B},\CC)$ is an isomorphism as vector spaces. In fact, if $\sigma\gamma=0$ for a $q$-cochain $\g$, then, by \eqref{sigma},
$\g_{\la_1,\cdots,\la_q}(a_1,\cdots,a_q)$ as a polynomial in $\la_q$ has a root $\la_q=-\mbox{$\sum_{i=1}^{q-1}$}\la_i$, namely, it is divided by $\mbox{$\sum_{i=1}^{q}$}\la_i$. Thus we get a $q$-cochain
\begin{eqnarray}
\gamma'_{\la_1,\cdots,\la_q}(a_1,\cdots,a_q)=(\mbox{$\sum_{i=1}^{q}$}\la_i)^{-1}\g_{\la_1,\cdots,\la_q}(a_1,\cdots,a_q),
\end{eqnarray}
and $\gamma=(\mbox{$\sum_{i=1}^{q}$}\la_i)\g'\in \pa {\t C}^q(\mathcal {B},\CC)$, which proves that $\sigma$ is injective. Hence the claim is true.

In the following we can identify $C^q(\mathcal {B},\CC)$ with $C'^q(\mathcal{B},\CC)$. We still call an element in $C'^q(\mathcal{B},\CC)$ a reduced $q$-cochain. By defining the operator $d':C'^q(\mathcal{B},\CC)\rightarrow C'^{q+1}(\mathcal{B},\CC)$ by $d' (\sigma\gamma)=\sigma d\gamma$, we have similar notions of reduced $q$-cocycle and $q$-coboundary. For convenience, we will abbreviate  $\g_{\la_1,\cdots,\la_q}(a_1,\cdots,a_q)|_{\la_q=-\la_1-\cdots-\la_{q-1}}$ to $\g_{\la_1,\cdots,\la_{q-1}}(a_1,\cdots,a_q)$.

Let $\psi '=\sigma \psi\in C'^2(\B,\CC)$ be a reduced 2-cochain. By \eqref{d} and \eqref{sigma},
\begin{eqnarray}\label{dd4}
&&(d \psi ')_{\la_1,\la_2}(a_1,a_2,a_3)\nonumber\\&&\ \ \ \ \ \ \ \, =-\psi '_{\la_1+\la_2}([a_{1\,\la_1}a_2],a_3)+\psi '_{-\la_2}([a_{1\,\la_1}a_3],a_2)-\psi '_{-\la_1}([a_{2\,\la_2}a_3],a_1),
\end{eqnarray}
for $a_1, a_2, a_3\in\B$. Define a reduce 1-cochain $f'=\sigma f\in C'^1(\B,\CC)$ by
\begin{eqnarray}\label{dd1}
f'(J_i)=(i+2)^{-1}\frac{d}{d\la}\psi'_{\la}(J_i,J_0)|_{\la=0}, \ \mbox{for}\ i\geq0.
\end{eqnarray}
Note that $f'(a)=f_{\la}(a)|_{\la=0}$. Thus $f'$ is simply a linear function from $\B$ to $\CC$, satisfying $f'(\pa a)=f_{\la}(\pa a)|_{\la=0}=-\la f_{\la}(a)|_{\la=0}=0$, and
\begin{eqnarray}\label{dd2}
(df')_{\la}(a_1,a_2)=-f'([a_{1\,\la}a_2]), \ \mbox{for}\ a_1,\ a_2\in\B.
\end{eqnarray}
If $\psi'$ is a reduced 2-cocycle, then $\g'=\psi'+d f'$ is a reduced 2-cocycle, equivalent to $\psi'.$ By \eqref{dd1}
and \eqref{dd2},
\begin{eqnarray}\label{dd3}
\frac{d}{d\la}\g'_{\la}(J_i,J_0)|_{\la=0}=0, \ \mbox{for}\ i\geq 0.
\end{eqnarray}
This, along with \eqref{dd4} and \eqref{skey-sym}, gives
\begin{eqnarray}
0&=&\frac{\pa}{\pa\la}(d\g')_{\la_1,\la}(J_i,J_k,J_0)|_{\la=-\la_1}\nonumber\\
&=&\frac{\pa}{\pa\la}\big(-\g'_{\la_1+\la}([J_{i\,{\la_1}}J_k],J_0)+\g'_{-\la}([J_{i\,{\la_1
}}J_0],J_k)-\g'_{-\la_1}([J_{k\,{\la}} J_0],J_i)\big)|_{\la=-\la_1}\nonumber\\
&=&\frac{\pa}{\pa\la}\big(((i+1)(\la+\la_1)+\la_1)\g'_{-\la}(J_i,J_k)-((k+1)(\la+\la_1)+\la)\g'_{-\la_1}(J_k,J_i)\big)|_{\la=-\la_1}.
\nonumber\\&=&(i+k+3)\g'_{\la_1}(J_i,J_k)-\la_1\frac{\pa}{\pa\la_1}\g'_{\la_1}(J_i,J_k).
\end{eqnarray}
Thus,
\begin{eqnarray}\label{2cocycle1}
 \g'_{\la}(J_i,J_k)=c_{i,k}\la^{i+k+3}, \ \mbox{for some}\ c_{i,k}\in\CC.
\end{eqnarray}
By \eqref{dd4} with $\g'$ in place of $\psi'$ and \eqref{2cocycle1},
\begin{eqnarray}\label{2co}
0&=& -\g '_{\la_1+\la_2}([J_{i\,\la_1}J_j],J_k)+\g '_{-\la_2}([J_{i\,\la_1}J_k],J_j)-\g '_{-\la_1}([J_{j\,\la_2}J_k],J_i)\nonumber\\
&=&-((j+1)\la_1-(i+1)\la_2)c_{i+j,k}(\la_1+\la_2)^{i+j+k+3}\nonumber\\&& +((i+1)\la_2+(i+k+2)\la_1)c_{i+k,j}(-\la_2)^{i+j+k+3}\nonumber\\&&
-((j+1)\la_1+(j+k+2)\la_2)c_{j+k,i}(-\la_1)^{i+j+k+3}.\label{2cocycle2}
\end{eqnarray}
Taking $i=j=0$ in \eqref{2cocycle2}, we get
\begin{eqnarray}\label{2cocycle3}
c_{0,k}(\la_1-\la_2)(\la_1+\la_2)^{k+3}=c_{k,0}\big((\la_2+(k+2)\la_1)(-\la_2)^{k+3}
-(\la_1+(k+2)\la_2)(-\la_1)^{k+3}\big).
\end{eqnarray}
Setting $\la_2=0$ in \eqref{2cocycle3} gives $c_{0,k}=(-1)^k c_{k,0}$.  Comparing coefficients of $\la_1^2\la_2^{k+2}$ in \eqref{2cocycle3} gives $c_{0,k}=c_{k,0}=0$ for $k\geq 1$. Setting $i=0$ in \eqref{2cocycle2}
and comparing coefficients of $\la_1^{j+k+4}$, we obtain $c_{j,k}=0$ for $j, k\geq 1$.
Thus, by \eqref{2cocycle1}, there exists a nonzero complex number $c$, such that
\begin{eqnarray}\label{2cocycle4}
 \g'_{\la}(J_0,J_0)=c\la^{3},  \ \ \g'_{\la}(J_i,J_j)=0, \ \mbox{for} \ i,\ j\geq 1.
\end{eqnarray}
Therefore, ${\rm dim \,H}^2(\B,\CC)=1$, which proves (1).

(2) Define an operator $\tau:\t C^q(\B,\CC_a)\rightarrow
\t C^{q-1}(\B,\CC_a)$ by
\begin{eqnarray}\label{7-3}
(\tau
\g)_{\la_1,\cdots,\la_{q-1}}(a_1,\cdots,a_{q-1})=(-1)^{q-1}\g_{\la_1,\cdots,\la_{q-1},\la}(a _1,\cdots,a_{q-1},J_0)|_{\la=0},
\end{eqnarray}
for $a_1,\cdots,a_{q-1}\in\B$. By the fact that $\pa\t C^q(\B,\CC_a)=(a+\sum_{i=1}^q\la_i)\t
C^q(\B,\CC_a)$ and \eqref{7-3},
\begin{eqnarray}\label{7-4}
((d\tau+\tau d)
\g)_{\la_1,\cdots,\la_{q}}(J_{n_1},\cdots,J_{n_q})&=&\big(\mbox{$\sum_{i=1}^q$}\la_i\big)\g_{\la_1,\cdots,\la_{q}}(J_{n_1},\cdots,J_{n_q})\nonumber\\
&\equiv& -a \g_{\la_1,\cdots,\la_{q}}(J_{n_1},\cdots,J_{n_q})\ (\mbox{mod}\
\pa\t C^q(\B,\CC_a).\ \
\end{eqnarray}
Suppose that $\g$ is
a $q$-cochain such that $d\g\in \pa\t C^{q+1}(\B,\CC_a)$, namely,
there is a $(q+1)$-cochain $\phi$ such that
$d\g=(a+\sum_{i=1}^{q+1}\la_i)\phi$. By \eqref{7-3}, $\tau d\g=(a+\sum_{i=1}^q\la_i)\tau\phi\in\pa\t
C^q(\B,\CC_a)$. By\eqref{7-4}, $\g\equiv-d(a^{-1}\tau\g)$ is a reduced coboundary because $a\neq 0$. Thus $ {\rm H}^q(\mathcal{B},\CC_a)=0$ for $q\geq 0$.
 This proves (2).

(3) Note that $\pa\t
C^q(\B,M_{\Delta,\a})=(\pa+\sum_{i=1}^q\la_i)\t C^q(\B,
M_{\Delta,\a})$. Similarly to the  proof of (2), we define an
operator $\kappa:C^q(\B,M_{\Delta,\a})\rightarrow
C^{q-1}(\B,M_{\Delta,\a})$ by
\begin{eqnarray*}\label{5-3}
(\kappa\g)_{\la_1,\cdots,\la_{q-1}}(a_1,\cdots,a_{q-1})=(-1)^{q-1}\g_{\la_1,\cdots,\la_{q-1},\la}(a_1,\cdots,a_{q-1},J_0)|_{\la=0},
\end{eqnarray*}
for $a_1,\cdots,a_{q-1}\in\B$. By Theorem \ref{theo},
\begin{eqnarray}\label{5-4}
&&((d\kappa+\kappa d)
\g)_{\la_1,\cdots,\la_{q}}(J_{n_1},\cdots,J_{n_q})\nonumber\\&&\ \ \ =J_{0\,\la}\g_{\la_1,\cdots,\la_{q}}(J_{n_1},\cdots,J_{n_q})|_{\la=0}+
\big(\mbox{$\sum\limits_{i=1}^q$}\la_i\big)\g_{\la_1,\cdots,\la_{q}}(J_{n_1},\cdots,J_{n_q})\nonumber\\&&\ \ \ = \big(\pa+\a+\mbox{$\sum\limits_{i=1}^q$}\la_i\big)\g_{\la_1,\cdots,\la_{q}}(J_{n_1},\cdots,J_{n_q})\nonumber\\&&\ \ \
\equiv \a \g_{\la_1,\cdots,\la_{q}}(J_{n_1},\cdots,J_{n_q}) \ (\mbox{mod}\
\pa\t C^q(\B, M_{\Delta,\a})).
\end{eqnarray}
If $\g$ is a reduced $q$-cocycle, then there is a $(q+1)$-cochain $\varphi$ such that $d\g=\pa \varphi=(\pa+\sum_{i=1}^{q+1}\la_i)\varphi$. In this case, $\kappa d\g=(\pa+\sum_{i=1}^q\la_i)\kappa\phi\in\pa\t
C^q(\B,M_{\Delta,\a})$. It follows from \eqref{5-4}
that $\g\equiv d(\a^{-1}\kappa\g)$ is a reduced $q$-coboundary, since we assume $\a\neq
0$. Hence ${\rm H}^q(\mathcal{B},M_{\Delta,\a})=0$ for $q\geq 0$.

This completes the proof of Theorem \ref{thm-5}. \QED

\begin{cor} There is a unique nontrivial universal central extension $\t\B=\B\oplus \CC {\mathfrak c}$ of the Lie conformal algebra $\B$, satisfying
\begin{eqnarray*}
[J_{0\,\la} J_0]&=&(\pa+2\la)J_{0}+\la^3 {\mathfrak c},\\
{[J_{i\,\la} J_j]}&=&((i+1)\pa+(i+j+2)\la)J_{i+j}, \ \mbox{for}\ i, j> 0.
\end{eqnarray*}
\end{cor}

\begin{rem}\rm The formal distribution Lie algebra corresponding to $\t\B$ is a well-known Lie algebra of Block type studied in \cite{SXX2,WT}.
\end{rem}

\section{Vertex Poisson algebra structure associated to $\B$}

Denote $V=\mbox{$\bigoplus_{i\in\Z^+}\CC J_i$}$. Thus the Lie conformal algebra of Block type
$\B$ is a free $\CC[\pa]$-module over $V$.
By \eqref{j-product}, the $\la$-bracket \eqref{bracket} is equivalent to the following  $j$-products
\begin{eqnarray}\label{jpts}
&&J_{i\,{(0)}}J_k=(i+1)\partial J_{i+k},\ \
J_{i\,{(1)}}J_k= (i+k+2)J_{i+k}, \ \ J_{i\,{(n)}}J_k=0,
\end{eqnarray}
for $i,k\in\Z^+$, $n\geq 2.$ Define a linear map $Y_-(\cdot, z)$ from $V$ to $z^{-1}({\rm End\,}
V)[[z^{-1}]]$ by
\begin{eqnarray}\label{y}
Y_-(J_i,z)=\mbox{$\sum\limits_{n\in\Z^+}$}J_{i\,{(n)}}z^{-n-1}, \ \mbox{for}\ i\geq 0.
\end{eqnarray}
From \eqref{jpts} and \eqref{y}, we have
\begin{eqnarray}\label{quadr}
Y_-(J_i,z)J_k=(i+1)\partial J_{i+k}z^{-1}+(i+k+2)J_{i+k}z^{-2}, \ \mbox{for} \ i,k\in\Z^+.
\end{eqnarray}
Extending the map $Y_-(\cdot, z)$ to the whole $\B=\CC[\pa]\otimes_\CC V$ by
\begin{eqnarray}\label{quadr3}
Y_-(f(\pa)J_i,z)\big(\pa^m
J_k\big)=f(d/dz)\mbox{$\sum\limits_{l=0}^{m}$}(-1)^l\pa^{m-l}(d/dz)^l
Y_-(J_i,z)J_k,
\end{eqnarray}
then $(\mathcal{B}, Y_-,\pa)$ forms a free vertex Lie algebra. Indeed, \eqref{quadr} and \eqref{quadr3} guarantee axioms \eqref{vHL0+} and \eqref{vHL1+}. It suffices to check \eqref{vHL2+} and \eqref{vHL55} on the generators. For $i,j\in\Z^+$,
\begin{eqnarray}\label{**}
{\rm Sing}(e^{z\pa}Y_-(J_j,-z)J_i)&=&(j+1)\partial J_{i+j}(-z)^{-1}+(i+j+2)J_{i+j}z^{-2}+(i+j+2)\pa J_{i+j}z^{-1}\nonumber\\
&=&(i+1)\partial J_{i+j}z^{-1}+(i+j+2)J_{i+j}z^{-2} \nonumber\\&=&Y_-(J_i,z)J_j,
\end{eqnarray}
which proves \eqref{vHL2+}. Furthermore, for any $i,j,k \in\Z^+$,
\begin{eqnarray}\label{v1}
&& Y_-(J_i,z_1) Y_-(J_j,z_2)J_k\nonumber\\
&&\ \  =(j+1)(i+1)\pa^2 J_{i+j+k}z_1^{-1}z_2^{-1}+(j+1)(2i+j+k+3)\pa J_{i+j+k}z_1^{-2}z_2^{-1}\nonumber\\
&&\ \ \ \ \ +2(j+1)(i+j+k+2)J_{i+j+k}z_1^{-3}z_2^{-1}+(i+1)(j+k+2)\pa J_{i+j+k} z_1^{-1}z_2^{-2}\nonumber\\&&\ \ \ \ \ +(j+k+2)(i+j+k+2)J_{i+j+k}z_1^{-2}z_2^{-2},
\end{eqnarray}
and
\begin{eqnarray}\label{***}
&&{\rm Sing}\big(\mbox{$\sum_{n\geq 0}$}(z_1-z_2)^{-n-1}Y_-(J_{i\,{(n)}}J_j,z_{2})\big)J_k\nonumber\\
&&\ ={\rm Sing}\big((z_1-z_2)^{-1}Y_-(J_{i\,{(0)}}J_j,z_{2})+(z_1-z_2)^{-2}Y_-(J_{\,i{(1)}}J_j,z_{2})\big)J_k\nonumber\\
&&\  ={\rm Sing}\big((i+1)(z_1-z_2)^{-1}Y_-(\pa J_{i+j},z_{2})+(i+j+2)(z_1-z_2)^{-2}Y_-( J_{i+j},z_{2})\big)J_k\nonumber\\
&&\ =-(i+1)(i+j+1)\pa J_{i+j+k}(z_1^{-1}z_2^{-2}+z_1^{-2}z_2^{-1})+(i+j+2)(i+j+1)\pa J_{i+j+k}z_1^{-1}z_2^{-2}\nonumber\\
&&\ \ \ \  -2(i+1)(i+j+k+2)J_{i+j+k}(z_1^{-1}z_2^{-3}+z_1^{-2}z_2^{-2}+z_1^{-3}z_2^{-1})\nonumber\\
&&\ \ \ \ +(i+j+2)(i+j+k+2)J_{i+j+k}(z_1^{-2}z_2^{-2}+2z_1^{-3}z_2^{-1})\nonumber\\
&&\ =Y_-(J_i,z_1) Y_-(J_j,z_2)J_k-Y_-(J_j,z_1) Y_-(J_i,z_2)J_k,
\end{eqnarray}
where the last equality follows from \eqref{v1}. This proves \eqref{vHL55} and thus $(\mathcal{B}, Y_-,\pa)$ is a vertex Lie algebra, called the {\it vertex Lie algebra of Block type}.
\begin{rem}{\rm By \cite{GD}, the vertex Lie algebra $(\mathcal{B}, Y_-,\pa)$ is equivalent to a Novikov algebra $(V,\circ)$ with the operation $\c$ defined by
\begin{eqnarray}\label{GB2}
J_i\circ J_j=(j+1)J_{i+j}, \ \mbox{for}\
\ i,j\geq 0.
\end{eqnarray}
Moreover, the commutator on $(V,\circ)$ of the form
\begin{eqnarray}
 [J_i,J_j]=(j-i)J_{i+j},\ \ \mbox{for}\
\ i,j\geq 0,
\end{eqnarray}
make $V$ into a Lie algebra. The Lie algebra $(V,[\cdot,\cdot])$ is a subalgebra of the one-sided Witt algebra $W_1^+$ with a basis $\{J_i|i\geq -1\}$, which occurs in the study of conformal field theory.}
\end{rem}

 By a {\it differential algebra}, we mean a commutative associative algebra $A$ (with $1$) equipped with a derivation $\pa$, denoted by $(A,\pa)$. A subset $U$ of $A$ is said to generate $A$ as a differential algebra if $\pa^n U$ for $n\geq 0$ generate $A$ as an algebra.

The following notion of vertex Poisson algebra is due to \cite{FB}.
\begin{defi} {\rm A vertex Poisson algebra is a differential algebra $(A,\pa)$ equipped with a vertex Lie algebra structure $(Y_-,\pa)$ such that for $a,b,c\in A$,
\begin{eqnarray}\label{vpa}
Y_-(a,z)(bc)=(Y_-(a,z)b)c+b(Y_-(a,z)c),
\end{eqnarray}
where $Y_-(a,z)=\mbox{$\sum_{n\geq 0}$}a_n z^{-n-1}$.}
\end{defi}

In terms of components, relation \eqref{vpa} is equivalent to
\begin{eqnarray}\label{vpa+}
a_n(bc)=(a_nb)c+b(a_nc) \ \mbox{for}\ a,b,c\in A,\ n\geq 0,
\end{eqnarray}
namely, these $a_n$ are derivations of $A$. Therefore,
\begin{eqnarray}\label{vpa++}
Y_-(a,z)\in z^{-1}({\rm\,Der} A)[[z^{-1}]],\ \ \mbox{for}\ a\in A,
\end{eqnarray}
which implies $Y_-(a,z)1=0$. Then $Y_-(1,z)a=0$ by \eqref{vHL2+}, namely, $Y_-(1,z)=0.$

 Let $U$ be a vector space. Denote by
$A=S(\CC[\pa]\otimes U)$
the symmetric algebra over the free $\CC[\pa]$-module $\CC[\pa]\otimes U$. The operator $\pa$ can be uniquely extended to a derivation of $A$. Then $(A,\pa)$ forms a differential algebra, which is referred to as {\it the free differential algebra} over $U$. Following \cite{Li1}, a {\it week pre-vertex Poisson structure} on $A$ is a linear map $Y_-^0$ from $U\times U$ to $z^{-1}A[ z^{-1}]$ such that
\begin{eqnarray}\label{fda+}
Y_-^0(a,z)b={\rm Sing}(e^{z\pa}Y_-^0(b,-z)a), \ \mbox{for}\ a,b\in A.
\end{eqnarray}
By \cite[Proposition 3.10]{Li1}, the operator $Y_-^0$ can be uniquely extended to a linear map
\begin{eqnarray*}
Y_-:A\rightarrow {\rm Hom}(A, z^{-1}A[ z^{-1}]),\ \ \
a\mapsto Y_-(a,z)=\mbox{$\sum\limits_{n\in\Z^+}$}a_n z^{-n-1}
\end{eqnarray*} such that for $a,b\in A$,
\begin{eqnarray}
Y_-(a,z)&=&\mbox{$\sum_{n\geq 0}$}a_n z^{-n-1}\in z^{-1}({\rm Der} A)[[z^{-1}]],\\
{[\pa, Y_-(a,z)]}&=&{Y_-(\pa a,z)}=\frac{d}{dz}Y_-(a,z),\\
{Y_-(a,z)b}&=& {\rm Sing} \big(e^{z\pa}Y_-(b,-z)a\big).
\end{eqnarray}

The following result, due to \cite[Theorem 3.11]{Li1}, gives a general construction of vertex Poisson algebras from free differential algebras.
\begin{thm} \label{theorem} Let $A$ be the free differential algebra over $U$ and let $Y_-^0$ be a week pre-vertex Poisson structure on $A$ such that for $u,v,w\in U$ from an ordered basis of $U$ with $u\leq v\leq w$,
\begin{eqnarray}\label{th1}
\t Y_-^0(u,z_1) Y_-^0(v,z_2)w-\t Y_-^0(v,z_2) Y^0(u,z_1)w
={\rm Sing}\big(e^{z_2\pa}\t Y_-^0(w,-z_{2}) Y_-^0(u,z_1-z_2)v\big),
\end{eqnarray}
where $\t Y_-^0(u,z)\in z^{-1}({\rm Der})[[z^{-1}]]$ is uniquely determined by
\begin{eqnarray}\label{th11}
\t Y_-^0(u,z) e^{z_1\pa}v=e^{z_1\pa}e^{-z_1\frac d{dz}}Y_-^0(u,z)v, \ \mbox{for}\ {v\in U}.
\end{eqnarray}
Then $Y_-^0$ uniquely extends to a vertex Poisson structure $Y_-$ on $A$.
\end{thm}

The following result is an application of the above theorem.
\begin{pro} Let $(\B,Y_-,\pa)$ be the vertex Lie algebra of Block type and $A=S(\B)$ the free differential algebra over $V=\mbox{$\bigoplus\limits_{i\in\Z^+}$}\CC J_i$. Define
\begin{eqnarray}\label{vpa-b}
Y_-^0(J_i,z)J_j=(i+1)\partial J_{i+j}z^{-1}+(i+j+2)J_{i+j}z^{-2}, \ \mbox{for}\ i,j\geq 0.
\end{eqnarray}
Then $Y_-^0$ uniquely extends to a vertex Poisson algebra structure on $A$.
\end{pro}
\pr ~ By \eqref{quadr} and \eqref{vpa-b}, $Y_-^{0}(J^i,z)J^j= Y_-(J^i,z)J^j$ for $i,j\in\Z^+$. Thus $Y_-^0$ satisfies \eqref{fda+} and so it is a week pre-vertex Poisson structure on $A$. By \eqref{th11} and \eqref{vpa-b}, for $i, j, k\in\Z^+$,
\begin{eqnarray}\label{vpa-b++}
&&\t Y_-^0(J_i,z_1) Y_-^0(J_j,z_2)J_k \nonumber\\
&&\ \ \ \ \ \ =(j+1)(i+1)\pa^2 J_{i+j+k}z_1^{-1}z_2^{-1}+(j+1)(2i+j+k+3)\pa J_{i+j+k}z_1^{-2}z_2^{-1}\nonumber\\
&&\ \ \ \ \ \ \ \ \ \,+2(j+1)(i+j+k+2)J_{i+j+k}z_1^{-3}z_2^{-1}+(i+1)(j+k+2)\pa J_{i+j+k} z_1^{-1}z_2^{-2}\nonumber\\&&\ \ \ \ \ \ \ \ \ \, +(j+k+2)(i+j+k+2)J_{i+j+k}z_1^{-2}z_2^{-2}.
\end{eqnarray}
By \eqref{th11}--\eqref{vpa-b++},
\begin{eqnarray*}
&&{\rm Sing}\big(e^{z_2\pa}\t Y_-^0(J_k,-z_{2}) Y_-^0(J_i,z_1-z_2)J_j\big)\\
&=&{\rm Sing}\big(e^{z_2\pa}Y_-^0(J_k,-z_{2})\big((i+1)\pa J_{i+j}(z_1-z_2)^{-1}+(i+j+2)J_{i+j}(z_1-z_2)^{-2}\big)\big)\\
&=&{\rm Sing}\big((i+1)Y_-(\pa J_{i+j},z_2)J_k(z_1-z_2)^{-1}+(i+j+2)Y_-( J_{i+j},z_2)J_k(z_1-z_2)^{-2}\big)\\
&=&{\rm Sing}\big(-(i+1)\big((i+j+1)\pa J_{i+j+k}z_2^{-2}+2(i+j+k+2)z_2^{-3}J_{i+j+k}\big)(z_1-z_2)^{-1}\\
&&+(i+j+2)\big((i+j+1)\pa J_{i+j+k}z_2^{-1}+(i+j+k+2)z_2^{-2}J_{i+j+k}\big)(z_1-z_2)^{-2}\big)\\
&=& -(i+1)(i+j+1)\pa J_{i+j+k}(z_2^{-2}z_1^{-1}+z_1^{-2}z_{2}^{-1})+(i+j+2)(i+j+1)\pa J_{i+j+k} z_2^{-1}z_1^{-2}\\&&
-2(i+1)(i+j+k+2)J_{i+j+k}(z_2^{-3}z_1^{-1}+z_2^{-2}z_1^{-2}+z_2^{-1}z_1^{-3})\\
&& +(i+j+2)(i+j+k+2)J_{i+j+k}(z_2^{-2}z_1^{-2}+2z_1^{-3}z_2^{-1})\\
&=&\t Y_-^0(L^i,z_1) Y_-^0(L^j,z_2)L^k-\t Y_-^0(L^j,z_2) Y_-^0(L^i,z_1)L^k.
\end{eqnarray*}
This proves the result by Theorem \ref{theorem}. \QED

 \vskip10pt
\small\noindent{\bf Acknowledgements~}~{\footnotesize
This work was supported by National Natural Science
Foundation grants of China (11301109) and
the Research Fund for the Doctoral Program of Higher Education
(20132302120042).}

\end{document}